\documentclass[11pt]{article}
\usepackage{geometry}                
\usepackage{pgf,tikz}
\usetikzlibrary{arrows}
\usepackage{graphicx}
\usepackage{amssymb,amsbsy,amsmath,amscd,amsthm}
\usepackage{epstopdf}
\usepackage{tikz}
\usetikzlibrary{arrows}

\usepackage[applemac]{inputenc}

\usepackage[english,francais]{babel}

\usepackage{graphicx}

\newtheorem{theorem}{Theorem}
\newtheorem{remarque}{Remark}

\newtheorem{lemma}{Lemma}

\newtheorem{corollary}{Corollary}

\newcommand{\R}{\mathbb R}
\newcommand{\PP}{\mathbb P}
\newcommand{\E}{\mathbb E}

\definecolor{ffqqqq}{rgb}{1,0,0}
\definecolor{qqwuqq}{rgb}{0,0.39,0}
\definecolor{xdxdff}{rgb}{0.49,0.49,1}
\definecolor{uququq}{rgb}{0.25,0.25,0.25}
\definecolor{zzttqq}{rgb}{0.6,0.2,0}
\definecolor{cqcqcq}{rgb}{0.75,0.75,0.75}

\title{\textbf{Adaptive estimation of convex and polytopal density support}}

\author{Victor-Emmanuel Brunel\\
   CREST-ENSAE, Malakoff, France / University of Haifa, Israel\\
  \texttt{victor.emmanuel.brunel@ensae-paristech.fr}}

\date{}                                           

\begin{document}

\maketitle{}

\selectlanguage{english}
\begin{abstract}
    We estimate the support of a uniform density, when it is assumed to be a convex polytope or, more generally, a convex body in $\R^d$. In the polytopal case, we construct an estimator achieving a rate which does not depend on the dimension $d$, unlike the other estimators that have been proposed so far. For $d\geq 3$, our estimator has a better risk than the previous ones, and it is nearly minimax, up to a logarithmic factor. We also propose an estimator which is adaptive with respect to the structure of the boundary of the unknown support.

\end{abstract}

\section{Introduction}\label{sec:1}
Assume we observe a sample of $n$ i.i.d. random variables $X_i, i=1,\ldots, n$, with uniform distribution on some subset $G$ of $\R^d, d\geq 2$. We are interested in the problem of estimation of $G$. In particular, this problem is of interest in detection of abnormal behavior, cf. Devroye and Wise \cite{DevroyeWise1980}. In image recovering, when an object is only partially observed, e.g. if only some pixels are available, one would like to recover the object as accurately as possible. When $G$ is known to be compact and convex, the convex hull of the sample is quite a natural estimator. The properties of this random subset of $\R^d$ have been extensively studied since the early 1960's, from a geometric and probabilistic prospective. The very original question associated to this object was the famous Sylvester four-point problem: what is the probability that one of the four points chosen at random in the plane is inside the triangle formed by the three others ? We refer to \cite{BaranySyl} for a historical survey and extensions of Sylvester problem. Of course, this question is not well posed, since the answer should depend on the probability measure of those four points, and the many answers that were proposed, in the late 18th century, accompanied the birth of a new field: stochastic geometry. Rényi and Sulanke \cite{26,27} studied some basic properties of the convex hull of $X_i, i=1,\ldots, n$ when $G$ is a compact and convex subset of the plane ($d=2$). More specifically, if this convex hull is denoted by $CH(X_1,\ldots,X_n)$, its number of vertices by $V_n$ and its missing area $|G\backslash CH(X_1,\ldots,X_n)|$ by $A_n$, they investigated the asymptotics of the expectations $\E[V_n]$ and $\E[A_n]$. Their results are highly dependent on the structure of the boundary of $G$. The expected number of vertices is of the order $n^{1/3}$ when the boundary of $G$ is smooth enough, and $r\ln n$ when $G$ is a convex polygon with $r$ vertices, $r\geq 3$. The expected missing area is of the order $n^{-2/3}$ in the first case and, if $G$ is a square, it is of the order $(\ln n)/n$. May the square be arbitrarily large or small, only the constants and not the rates are affected, by a scale factor. Rényi and Sulanke \cite{26,27} provided asymptotic evaluations of these expectations with the explicit constants, up to two or three terms. In 1965, Efron \cite{13} showed a very simple equality which connects the expected value of the number of vertices $V_{n+1}$ and that of the missing area $A_n$. Namely, if $|G|$ stands for the area of $G$, one has
\begin{equation}
	\label{EfronId} \E[A_n]=\frac{|G|\E[V_{n+1}]}{n+1},
\end{equation}
independently of the structure of the boundary of $G$. In particular, \eqref{EfronId} allows to extend the results of \cite{26,27} about the missing area to any convex polygon with $r$ vertices. If $G$ is such a polygon, $\E[A_n]$ is of the order $r(\ln n)/n$, up to a factor of the form $c|G|$, where $c$ is positive and does not depend on $r$ or $G$. More recently, many efforts were made to extend these results to dimensions 3 and more. We refer to \cite{Wieacker1987}, \cite{15bis}, \cite{Dwyer1988} and the references therein. Notably, Efron's identity \eqref{EfronId} holds in any dimension if $G\subseteq\R^d$ is a compact and convex set and $|G|$ is its Lebesgue measure.

Bàràny and Larman \cite{3'} (see \cite{3} for a review) proposed a generalization of these results with no assumption on the structure of the boundary of $G$. They considered the $\varepsilon$-wet part of $G$, denoted by $G(\varepsilon)$ and defined as the union of all the caps of $G$ of volume $\varepsilon|G|$, where a cap is the intersection of $G$ with a half space. Here, $0\leq\varepsilon\leq 1$. This notion, together with that of floating body (defined as $G\backslash G(\varepsilon)$) had been introduced by Dupin \cite{Dupin} and, later, by Blaschke \cite{Blaschke}. In \cite{3'}, the authors prove that the expected missing volume of the convex hull of independent random points uniformly distributed in a convex body is of the order of the volume of the $1/n$-wet part. Then the problem of computing this expected missing volume becomes analytical, and learning about its asymptotics reduces to analyzing the properties of the wet part, which have been studied extensively in convex analysis and geometry. In particular we refer to \cite{SchuttWerner}, \cite{Schutt94}, \cite{Werner06}, \cite{toappear1} and the references therein. In particular, it was shown that if the boundary of the convex body $G$ is smooth enough, then the expected missing volume is of the order $n^{-2/(d+1)}$, and if $G$ is a polytope, the order is $(\ln n)^{d-1}n^{-1}$.

All these works were developed in a geometric and probabilistic prospective. No efforts were made at this stage to understand whether the convex hull estimator is optimal if seen as an estimator of the set $G$. Only in the 1990's, this question was invoked in the statistical literature. Mammen and Tsybakov \cite{21'} showed that under some restrictions on the volume of $G$, the convex hull is optimal in a minimax sense (see the next section for details). Korostelev and Tsybakov \cite{KTlectureNotes1993} give a detailed account of the topic of set estimation. See also \cite{8'}, \cite{8''}, \cite{10}, \cite{Guntuboyina}, for an overview of recent developments about estimation of the support of a probability measure. A different model was studied in \cite{Brunel13}, where we considered estimation of the support of the regression function. We built an estimator which achieves a speed of convergence of the order $(\ln n)/n$ when the support is a polytope in $\R^d,d\geq 2$. Moreover, we proved that no estimator can achieve a better speed of convergence, so the logarithmic factor cannot be dropped. Although our estimator depends on the knowledge of the number of vertices $r$ of the true support of the regression function, we proposed an adaptive estimator, with respect to $r$, which achieves the same speed as for the case of known $r$.

However, to our knowledge, when one estimates the support of a uniform distribution, there are no results about optimality of the convex hull estimator when that support is a general convex set. In particular, when no assumptions on the location and on the structure of the boundary are made, it is not known if a convex set can be uniformly consistently estimated. In addition, the case of polytopes has not been investigated. Intuitively, the convex hull estimator can be improved, and the logarithmic factor can be, at least partially, dropped. Indeed, a polytope with a given number of vertices is completely determined by the coordinates of its vertices, and therefore belongs to some parametric family. This paper is organized as follows. In Section \ref{ssec:1.2} we give all notation and definitions. In Section \ref{sec:2} we propose a new estimator of $G$ when it is assumed to be a polytope and its number of vertices is known. We show that the risk of this estimator is better than that of the convex hull estimator, and achieves a rate independent of the dimension $d$. In Section \ref{sec:3}, we show that in the general case, if no other assumption than compactness, convexity and positive volume is made on $G$, then the convex hull estimator is optimal in a minimax sense. In Section \ref{sec:4} we construct an estimator which is adaptive to the shape of the boundary of $G$, i.e. which detects, in some sense, whether $G$ is a polytope or not and, if yes, correctly estimates its number of vertices. Section \ref{sec:6} is devoted to the proofs.

\section{Notation and Definitions}\label{ssec:1.2}

Let $d\geq 2$ a positive integer. Denote by $\rho$ the Euclidean distance in $\R^d$ and by $B_2^d$ the Euclidean unit closed ball in $\R^d$.

For brevity, we will a call convex body any compact and convex subset of $\R^d$ with positive volume, and we will call a polytope any compact convex polytope in $\R^d$ with positive volume. For an integer $r\geq d+1$, denote by $\mathcal P_r$ the class of all convex polytopes in $[0,1]^d$ with at most $r$ vertices. Denote also by $\mathcal K$ the class of all convex bodies in $\R^d$.

If $G$ is a closed subset of $\R^d$ and $\epsilon$ is a positive number, we denote by $G^\epsilon$ the set of all $x\in\R^d$ such that $\rho(x,G)\leq\epsilon$ or, in other terms, $G^\epsilon=G+\epsilon B_2^d$. If $G$ is any set, $I(\cdot\in G)$ stands for the indicator function of $G$.

The Lebesgue measure on $\R^d$ is denoted by $|\cdot|$ (for brevity, we do not indicate explicitly the dependence on $d$). If $G$ is a measurable subset of $\R^d$, we denote respectively by $\PP_G$ and $\E_G$ the probability measure of the uniform distribution on $G$ and the corresponding expectation operator, and we still use the same notation for the $n$-product of this distribution if there is no possible confusion. When necessary, we add the superscript $^{\otimes n}$ for the $n$ product. We will use the same notation for the corresponding outer probability and expectation when it is necessary, to avoid measurability issues. The Nikodym pseudo distance between two measurable subsets $G_1$ and $G_2$ of $\R^d$ is defined as the Lebesgue measure of their symmetric difference, namely $|G_1\triangle G_2|$.

A subset $\hat G_n$ of $\R^d$, whose construction depends on the sample is called a set estimator or, more simply, an estimator.

Given an estimator $\hat G_n$, we measure its accuracy  on a given class of sets in a minimax framework. The risk of $\hat G_n$ on a class $\mathcal C$ of Borel subsets of $\R^d$ is defined as
\begin{equation}
	\label{defriskest}\mathcal R_n(\hat G_n ; \mathcal C) = \sup_{G\in\mathcal C}\mathbb E_G[|G\triangle\hat G_n|].\tag{$*$}
\end{equation}
The rate (a sequence depending on $n$) of an estimator on a class $\mathcal C$ is the speed at which its risk converges to zero when the number $n$ of available observations tends to infinity. For all the estimators defined in the sequel we will be interested in upper bounds on their risk, in order to get information about their rate. For a given class of subsets $\mathcal C$, the minimax risk on this class when $n$ observations are available is defined as
\begin{equation}
\label{defriskminimax}\mathcal R_n(\mathcal C)=\inf_{\hat G_n} \mathcal R_n(\hat G_n ; \mathcal C), \tag{$**$}
\end{equation}
where the infimum is taken over all set estimators depending on $n$ observations. If $\mathcal R_n(\mathcal C)$ converges to zero, we call the minimax rate of convergence on the class $\mathcal C$ the speed at which $\mathcal R_n(\mathcal C)$ tends to zero.
For a given class $\mathcal C$ of subsets of $\R^d$, it is interesting to provide a lower bound for $\mathcal R_n(\mathcal C)$. By definition, no estimator can achieve a better rate on $\mathcal C$ than that of the lower bound. This bound gives also information on how close the risk of a given estimator is to the minimax risk. If the rate of the upper bound on the risk of an estimator matches the rate of the lower bound on the minimax risk on the class $\mathcal C$, then the estimator is said to have the minimax rate of convergence on this class.


For two quantities $A$ and $B$, and a parameter $\vartheta$, which may be multidimensional, we will write $A\lesssim_ \vartheta B$ (respectively $A\gtrsim_\vartheta B$) to say that for some constant positive constant $c(\vartheta)$ which depends on $\vartheta$ only, one has $A\leq c(\vartheta)B$ (respectively $A\geq c(\vartheta)B$). If we put no subscript under the signs $\lesssim$ or $\gtrsim$, this means that the involved constant is universal, i.e. depends on no parameter.


\section{Estimation of polytopes}\label{sec:2}
\subsection{Upper bound}\label{ssec:2.1}
Let $r\geq d+1$ be a known integer. Assume that the underlying set $G$, denoted by $P$ in this section, is in $\mathcal P_r$. The likelihood of the model, seen as a function of the compact set $G'\subseteq\R^d$, is defined as follows, provided that $G'$ has a positive Lebesgue measure:
\begin{equation*}
	L(X_1,\ldots,X_n,G')=\prod_{i=1}^n \frac{I(X_i\in G')}{|G'|}.
\end{equation*}
Therefore, maximization of the likelihood over a given class $\mathcal C$ of candidates, $\max_{G'\in\mathcal C}L(X_1,\ldots,X_n,G')$, is achieved when $G'$ is of minimum Lebesgue measure among all sets of $\mathcal C$ containing all the sample points. When $\mathcal C$ is the class of all convex subsets of $\R^d$, the maximum likelihood estimator is unique, and it is the convex hull of the sample. As we discussed above, this estimator has been extensively studied. In particular, using Efron's identity \eqref{EfronId}, it turns out that its expected number of vertices is of the order of $r(\ln n)^{d-1}$. However, the unknown polytope $P$ has no more than $r$ vertices. Hence, it seems reasonable to restrict the estimator of $P$ to have much less vertices; the class of all convex subsets of $\R^d$ is too large and we propose to maximize the likelihood over the smaller class $\mathcal P_r$.

Assume that there exists a polytope in $\mathcal P_n$ with the smallest volume among all polytopes of $\mathcal P_r$ containing all the sample points. Let $\hat P_n^{(r)}$ be such a polytope, i.e.
\begin{equation}
    \label{defest} \hat P_n^{(r)}\in\underset{P\in\mathcal P_r, X_i\in P, i=1,\ldots,n}{\operatorname{argmin}} \text{ } |P|.
\end{equation}
The existence of such a polytope is ensured by compactness arguments. Note that $\hat P_n^{(r)}$ needs not be unique.
The next theorem establishes an exponential deviation inequality for the estimator $\hat P_n^{(r)}$.
\begin{theorem}\label{Theorem1}
Let $r\geq d+1$ be an integer, and $n\geq 2$. Then,
\begin{equation*}
    \sup_{P\in\mathcal P_r}\mathbb P_{P}\left[n\left(|\hat P_n^{(r)}\triangle P|-\frac{4dr\ln n}{n}\right)\geq x\right] \lesssim_ d\mathrm e^{-x/2}, \forall x>0.
\end{equation*}
\end{theorem}

From the deviation inequality of Theorem~1 one can easily derive that the risk of the estimator $\hat P_n^{(r)}$ on the class $\mathcal P_r$ is of the order $\frac{\ln n}{n}$. Indeed, we have the next corollary.

\begin{corollary}\label{Corollary1} Let the assumptions of Theorem~1 be satisfied. Then, for any positive number $q$,
\begin{equation*}
\sup_{P\in\mathcal P_r} \mathbb E_{P}\left[|\hat P_n^{(r)}\triangle P|^q\right]\lesssim_ {d,q}\left(\frac{r\ln n}{n}\right)^q.
\end{equation*}
\end{corollary}
Corollary \ref{Corollary1} shows that the risk $\mathcal R_n(\hat P_n^{(r)} ; \mathcal P_r)$ of the estimator $\hat P_n^{(r)}$ on the class $\mathcal P_r$ is bounded from above by $\frac{r\ln n}{n}$, up to some positive constant which depends on $d$ only. Therefore we have the following upper bound for the minimax risk on the class $\mathcal P_r$:
\begin{equation}
	\label{ubound1} \mathcal R_n(\mathcal P_r) \lesssim_ d \frac{r\ln n}{n}.
\end{equation}
It is now natural to ask wether the rate $\frac{\ln n}{n}$ is minimax, i.e. whether it is possible to find a lower bound for $\mathcal R_n(\mathcal P_r)$ which converges to zero at the rate $\frac{\ln n}{n}$, or the logarithmic factor should be dropped. This question is discussed in the next subsection.

\subsection{The logarithmic factor}\label{ssec:5.2}

We conjecture that the logarithmic factor can be removed in the upper bound of $\mathcal R_n(\mathcal P_r), r\geq d+1$. Specifically, for the class of all convex polytopes with at most $r$ vertices, not necessarily included in the square $[0,1]^d$, which we denote by $\mathcal P_r^{all}$, we conjecture that, for the normalized version of the risk,
\begin{equation*}
	\mathcal Q_n(\mathcal P_r^{all})\lesssim_ d \frac{r}{n}.
\end{equation*}
What motivates our intuition is Efron's identity \eqref{EfronId}. Let us recall its proof, which is very easy, and instructive for our purposes.
Let the underlying set $G$ be a convex body in $\R^d$, denoted by $K$, and let $\hat K_n$ be the convex hull of the sample. Almost surely, $\hat K_n\subseteq K$, so
\begin{align}
	\E_K^{\otimes n}[|\hat K_n\triangle K|] & = \E_K^{\otimes n}[|K\backslash\hat K_n|] \nonumber \\
	& =  \E_K^{\otimes n}\left[\int_K I(x\notin\hat K_n)dx\right] \nonumber \\
	& = |K|\E_K^{\otimes n}\left[\frac{1}{|K|}\int_K I(x\notin\hat K_n)dx\right] \nonumber \\
	\label{EfronProof1} & = |K|\E_K^{\otimes n}\left[\PP_K[X\notin \hat K_n|X_1,\ldots,X_n]\right],
\end{align}
where $X$ is a random variable with the same distribution as $X_1$, and independent of the sample $X_1,\ldots,X_n$, and $\PP_K[\cdot|X_1,\ldots,X_n]$ denotes the conditional distribution given $X_1,\ldots,X_n$.
In what follows, we set $X_{n+1}=X$, so that we can consider the bigger sample $X_1,\ldots,X_{n+1}$. For $i=1,\ldots,n+1$, we denote by $\hat K^{-i}$ the convex hull of the sample $X_1,\ldots,X_{n+1}$ from which the $i$-th variable $X_i$ is withdrawn. Then $\hat K_n=\hat K^{-(n+1)}$, and by continuing \eqref{EfronProof1}, and by using the symmetry of the sample,
\begin{align}
	\E_K^{\otimes n}[|\hat K_n\triangle K|] & = |K|\PP_K^{\otimes n+1}[X_{n+1}\notin \hat K^{-(n+1)}] \nonumber \\
	& = \frac{|K|}{n+1}\sum_{i=1}^{n+1}\PP_K^{\otimes n+1}[X_i\notin \hat K^{-(i)}] \nonumber \\
	\label{EfronProof2} & =  \frac{|K|}{n+1}\sum_{i=1}^{n+1}\PP_K^{\otimes n+1}[X_i\in V(\hat K_{n+1})],
\end{align}
where $V(\hat K_{n+1})$ is the set of vertices of $\hat K_{n+1}=CH(X_1,\ldots,X_{n+1})$. Indeed, with probability one, the point $X_i$ is not in the convex hull of the $n$ other points if and only if it is a vertex of the convex hull of the whole sample. By rewriting the probability of an event as the expectation of its indicator function, one gets from \eqref{EfronProof2},
\begin{align*}
	\E_K^{\otimes n}[|\hat K_n\triangle K|] & = \frac{|K|}{n+1}\sum_{i=1}^{n+1}\E_K^{\otimes n+1}[I(X_i\in V(\hat K_{n+1}))] \nonumber \\
	& = \frac{|K|}{n+1}\E_K^{\otimes n+1}\left[\sum_{i=1}^{n+1}I(X_i\in V(\hat K_{n+1}))\right] \nonumber \\
	& =  \frac{|K|\E_K^{\otimes n+1}[V_{n+1}]}{n+1},
\end{align*}
where $V_{n+1}$ denotes the cardinality of $V(\hat K_{n+1})$, i.e. the number of vertices of the convex hull $\hat K_{n+1}$. Efron's equality is then proved.

It turns out that we can follow almost all the proof of this identity when the underlying set $G$ is a polytope, and when we consider the estimator developed in Section \ref{ssec:2.1}.
Let $r\geq d+1$ be an integer and $P\in\mathcal P_r^{all}$. Let $\hat P_n^{(r)}$ be the estimator defined in \eqref{defest}, where $\mathcal P_r$ is replaced by $\mathcal P_r^{all}$. In this section, we denote this estimator simply by $\hat P_n$. Note that this estimator does not satisfy the nice property $\hat P_n\subseteq P$, unlike the convex hull. However, by construction, $|\hat P_n|\leq |P|$, so $|P\triangle\hat P_n|\leq 2|P\backslash\hat P_n|$, and we have:
\begin{align}
	\E_P^{\otimes n}[|\hat P_n\triangle P|] & \leq 2\E_P^{\otimes n}[|P\backslash \hat P_n|] \nonumber \\
	& = 2|P|\E_P^{\otimes n}\left[\frac{1}{|P|}\int_P I(x\notin\hat P_n)dx\right] \nonumber \\
	\label{Efron'Proof1} & = 2|P|\E_P^{\otimes n}\left[\PP_P[X\notin \hat P_n|X_1,\ldots,X_n]\right],
\end{align}
where $X$ is a random variable of the same distribution as $X_1$, and independent of the sample $X_1,\ldots,X_n$, and $\PP_P[\cdot|X_1,\ldots,X_n]$ denotes the conditional distribution of $X$ given $X_1,\ldots,X_n$.
We set $X_{n+1}=X$, and we consider the bigger sample $X_1,\ldots,X_{n+1}$. For $i=1,\ldots,n+1$, we denote by $\hat P^{-i}$ the same estimator as $\hat P_n$, but this time based on the sample $X_1,\ldots,X_{n+1}$ from which the $i$-th variable $X_i$ is withdrawn. In other words, $\hat P^{-i}$ is a convex polytope with at most $r$ vertices, which contains the whole sample $X_1,\ldots,X_{n+1}$ but maybe the $i$-th variable, of minimum volume. Then, $\hat P_n=\hat P^{-(n+1)}$, and by continuing \eqref{Efron'Proof1},
\begin{align}
	\E_P^{\otimes n}[|\hat P_n\triangle P|] & \leq 2|P|\PP_P^{\otimes n+1}[X_{n+1}\notin \hat P^{-(n+1)}] \nonumber \\
	& = \frac{2|P|}{n+1}\sum_{i=1}^{n+1}\PP_P^{\otimes n+1}[X_i\notin \hat P^{-(i)}] \nonumber \\
	& = \frac{2|P|}{n+1}\E_P^{\otimes n+1}\left[\sum_{i=1}^{n+1}I(X_i\notin \hat P^{-(i)})\right] \nonumber \\
	\label{Efron'Proof2} & = \frac{2|P|\E_P^{\otimes n+1}[V'_{n+1}]}{n+1},
\end{align}
where $V'_{n+1}$ stands for the number of points $X_i$ falling outside of the polytope with at most $r$ vertices, of minimum volume, containing all the other $n$ points. Note that in this description we assume the uniqueness of such a polytope, which we conjecture to hold almost surely, as long as $n$ is large enough. It is not clear that if a point $X_i$ is not in $\hat P^{-i}$, then $X_i$ lies on the boundary of $\hat P_{n+1}$. However, if this was true, then almost surely $V'_{n+1}$ would be less or equal to $d+1$ times the number of facets of $\hat P_{n+1}$, since any facet is supported by an affine hyperplane of $\R^d$, which, with probability one, cannot contain more than $d+1$ points of the sample at a time. Besides, the maximal number of facets of a $d$ dimensional convex polytope with at most $r$ vertices is bounded by McMullen's upper bound \cite{McMullen}, \cite{Bronsted}, and we could have our conjecture proved. However, there might be some cases when some points $X_i$ are not in $\hat P^{-i}$, though they do not lay on the boundary of $\hat P_{n+1}$. So it may be of interest to work directly on the variable $V'_{n+1}$. This remains an open problem.

\subsection{Lower bound for the minimax risk in the case $d=2$}\label{ssec:5.3}

In the 2-dimensional case, we provide a lower bound of the order $1/n$, with a factor that is linear in the number of vertices $r$. Namely, the following theorem holds.

\begin{theorem}\label{theorem5}
Let $r\geq 10$ be an integer, and $n\geq r$. Assume $d=2$. Then,
\begin{equation*}
	\mathcal R_n(\mathcal P_r)\gtrsim \frac{r}{n}.
\end{equation*}
\end{theorem}

Combined with \eqref{ubound1}, this bound shows that, as a function of $r$, $\mathcal R_n(\mathcal P_r)$ behaves linearly in $r$ in dimension two. In greater dimensions, it is quite easy to show that $\mathcal R_n(\mathcal P_r)\gtrsim_{d} \frac{1}{n}$, but this lower bound does not show the dependency in $r$. However, the upper bound \eqref{ubound1} shows that $\mathcal R_n(\mathcal P_r)$ is at most linear in $r$.

\section{Estimation of convex bodies}\label{sec:3}

In this section we no longer assume that the unknown support $G$ belongs to a class $\mathcal P_r, r\geq d+1$, but only that it is a convex body and we write $G=K$. Denote by $\hat K_n$ the convex hull of the sample. The risk of this estimator cannot be bounded from above uniformly on the class $\mathcal K$, since by \eqref{EfronId} for any given $n$, $\E_K[|K\triangle\hat K_n|] \rightarrow \infty$ as $|K|\rightarrow \infty$. Moreover there is no uniformly consistent estimator on the class $\mathcal K$ of all convex bodies if the risk is defined by \eqref{defriskest}. The following result holds.
\begin{theorem}\label{theorem3}
	For all $n\geq 1$, the minimax risk \eqref{defriskminimax} on the class $\mathcal K$ is infinite:
	\begin{equation*}
		\mathcal R_n(\mathcal K)=+\infty.
	\end{equation*}
\end{theorem}
Therefore we will use another risk measure which is the normalized risk for an estimator $\tilde K_n$ of $K$, based of a sample of $n$ observations:
\begin{equation*}
	\mathcal Q_n(\tilde K_n ; \mathcal K) = \sup_{K\in\mathcal K}\mathbb E_K\left[\frac{|K\triangle\tilde K_n|}{|K|}\right].
\end{equation*}
Also define the normalized minimax risk on the class $\mathcal K$:
\begin{equation*}
	\mathcal Q_n(\mathcal K) = \inf_{\tilde K_n}\sup_{K\in\mathcal K}\mathbb E_K\left[\frac{|K\triangle\tilde K_n|}{|K|}\right],
\end{equation*}
where the infimum is taken over all estimators $\tilde K_n$ based on a sample of $n$ i.i.d. observations.
For the estimator $\hat K_n$ we do not provide a deviation inequality as in Theorem~1, but only an upper bound on the normalized risk.

\begin{theorem}\label{Theorem2}
Let $n\geq 2$ be an integer. Then,
\begin{equation*}
	\mathcal Q_n(\hat K_n ; \mathcal K) \lesssim_ d n^{-\frac{2}{d+1}}.
\end{equation*}
\end{theorem}
Note that this result gives a bound on $\E_K\left[\frac{|K\triangle\hat K_n|}{|K|}\right]$ that is uniform over \textit{all} convex bodies in $\R^d$, with no restriction on the location of the set $K$ (such as $K\subseteq[0,1]^d$) or on the volume of $K$, unlike in \cite{21'}. From Theorem \ref{Theorem2} and the lower bound of \cite{21'} (the lower bound of \cite{21'} is for the minimax risk, but the proof still holds for the normalized risk), we obtain the next corollary.

\begin{corollary}\label{corollary2}
	Let $n\geq 2$ be an integer. The normalized minimax risk on the class $\mathcal K$ satisfies
	\begin{equation*}
	n^{-\frac{2}{d+1}}\lesssim_ d \mathcal Q_n(\mathcal K) \lesssim_ d n^{-\frac{2}{d+1}},
	\end{equation*}
	and the convex hull has the minimax rate of convergence on $\mathcal K$, with respect to the normalized version of the risk.
\end{corollary}

Note that if instead of the class $\mathcal K$ we consider the class $\mathcal K_1$ of all convex bodies that are included in $[0,1]^d$, then $\forall K\in\mathcal K_1, |K|\leq 1$ and therefore, the risk of the convex hull estimator $\hat K_n$ on this class is bounded from above by $n^{-2/(d+1)}$, so
\begin{equation*}
	\mathcal R_n(\mathcal K_1) \lesssim_ d n^{-\frac{2}{d+1}}.
\end{equation*}
Besides, the lower bound that is given in \cite{21'} still holds for the class $\mathcal K_1$, and thus we have the following corollary.

\begin{corollary}\label{corollary3}
	Let $n\geq 2$ be an integer. The minimax risk on the class $\mathcal K_1$ satisfies
	\begin{equation*}
	n^{-\frac{2}{d+1}}\lesssim_ d \mathcal R_n(\mathcal K_1) \lesssim_ d n^{-\frac{2}{d+1}},
	\end{equation*}
	and the convex hull has the minimax rate of convergence on $\mathcal K_1$, with respect to the risk defined in \eqref{defriskminimax}.
\end{corollary}


\section{Adaptative estimation}\label{sec:4}

In Sections \ref{sec:2} and \ref{sec:3}, we proposed estimators which highly depend on the structure of the boundary of the unknown support. In particular, when the support was supposed to be polytopal with at most $r$ vertices, for some known integer $r$, our estimator was by construction also a polytope with at most $r$ vertices. Now we will construct an estimator which does not depend on any other knowledge than the convexity of the unknown support, and the fact that it is located in $[0,1]^d$. This estimator will achieve the same rate as the estimators of Section \ref{ssec:2.1} in the polytopal case, that is, $r\ln n/n$, where $r$ is the unknown number of vertices of the support, and the same rate, up to a logarithmic factor, as the convex hull which was studied in Section \ref{sec:3} in the case where the support is not polytopal, or is polytopal but with too many vertices. Note that if the support is a polytope with $r$ vertices, where $r$ is larger than $(\ln n)^{-1}n^\frac{d-1}{d+1}$, then the risk of the convex hull estimator $\hat K_n$ has a smaller rate than that of $\hat P_n^{(r)}$.
The idea which we develop here is the same as in \cite{Brunel13}, Theorem 6. The classes $\mathcal P_r, r\geq d+1$, are nested, that is, $\mathcal P_r\subseteq \mathcal P_{r'}$ as soon as $r\leq r'$. So it is better, in some sense, to overestimate the true number vertices of the unknown polytope $P$. Intuitively, it makes sense to fit some polytope with more vertices to $P$, while the opposite may be impossible (e.g. it is possible to fit a quadrilateral on any triangle, but not to fit a triangle on a square). We use this idea in order to select an estimator among the preliminary estimators $\hat P_n^{(r)}, r\geq d+1$, and $\hat K_n$.
Note that in \cite{Brunel13}, Theorem 6, the key tools for adaptation are the deviation inequalities for the preliminary estimators, but we do not have such an inequality for the "last" one, i.e. the convex hull $\hat K_n$. This induces a loss of precision in our estimation procedure. Namely, an extra logarithmic factor will appear.

Set $R_n=\lfloor n^{(d-1)/(d+1)}\rfloor$, where $\lfloor\cdot\rfloor$ stands for the integer part. Let $C$ be any positive constant greater than $16d+\frac{16}{d+1}$, and define
\begin{equation*}
	\hat r=\min\left\{ r\in\{d+1,\ldots,n\} : |\hat P_n^{(r)}\triangle\hat P_n^{(r')}|\leq \frac{Cr'\ln n}{n}, \forall r'=r,\ldots,n\right\}.
\end{equation*}
The integer $\hat r$ is well defined ; indeed, the set in the brackets in the last display is not empty, since the inequality is satisfied for $r=n$. As previously, denote by $V_n$ the number of vertices of $\hat K_n$, the convex hull of the sample. By definition of the convex hull, $\hat P_n^{(r)}=\hat K_n$, for all $r\geq V_n$. Therefore $\hat r\leq V_n$ almost surely.

The adaptive estimator is defined as follows.
\begin{equation*}
	\hat P_n^{adapt}= { \left\{
    \begin{array}{l}
        \hat P_n^{(\hat r)}   \mbox{ }\mbox{ if } \hat r\leq R_n   \vspace{3mm} \\
        \hat K_n, \mbox{ }\mbox{ otherwise}.
    \end{array}
	\right.}
\end{equation*}

Then, if we denote by $\mathcal P_{\infty}=\mathcal K_1$, we have the following theorem.

\begin{theorem}\label{theorem4}
Let $n\geq 2$. Let $\phi_{n,r}=\min\left(\frac{r\ln n}{n},(\ln n)n^{-\frac{2}{d+1}}\right)$, for all integers $r\geq d+1$ and $r=\infty$. Then,
\begin{equation*}
	\sup_{d+1\leq r\leq \infty}\sup_{P\in\mathcal P_r} \mathbb E_P\left[\phi_{n,r}^{-1}|\hat P_n^{adapt}\triangle P|\right] \lesssim_ d 1.
\end{equation*}
\end{theorem}

Thus, we show that one and the same estimator $\hat P_n^{adapt}$ attains the optimal rate, up to a logarithmic factor, simultaneously on all the classes $\mathcal P_r, d+1\leq r$ and on  the class $\mathcal K_1$ of all convex bodies in $[0,1]^d$.

\begin{remarque}
	The construction of $\hat r$ is inspired by Lepski's method \cite{21'0}. However, we cannot use the same techniques here since they need deviation inequalities for all the preliminary estimators, while we do not have such an inequality for $\hat K_n$. Our proof uses directly the properties of the convex hull estimator.
\end{remarque}

\section{Proofs}\label{sec:6}

\subsection{Proof of Theorem \ref{Theorem1}}

Let $r\geq d+1$ be an integer, and $n\geq 2$. Let $P_0\in\mathcal P_r$ and consider a sample $X_1,\ldots,X_n$ of i.i.d. random variables with uniform distribution on $P_0$.
For simplicity's sake, we will denote $\hat P_n$ instead of $\hat P_n^{(r)}$ in this Section.

Let $\hat P_n$ be the estimator defined in Theorem \ref{Theorem1}.
Let us define $\mathcal P_r^{(n)}$ as the class of all convex polytopes of $\mathcal P_r$ whose vertices lay on the grid $\left(\frac{1}{n}\mathbb Z\right)^d$, i.e. have as coordinates integer multiples of $1/n$.
We use the following lemma, whose proof can be found in \cite{Brunel13}.
\begin{lemma}
Let $r\geq d+1, n\geq 2$. There exists a positive constant $K_1$, which depends on $d$ only, such that for any convex polytope $P$ in $\mathcal P_r$ there is a convex polytope $P^*\in\mathcal P_{r}^{(n)}$ such that :
\begin{equation}
{\label{lem0} \left\{
    \begin{array}{l}
        |P^*\triangle P|\leq \frac{K_1}{n}     \vspace{3mm}\\
        P^*\subseteq P^{\sqrt d/n}, \mbox{ }\mbox{ }P\subseteq (P^*)^{\sqrt d/n}.
    \end{array}
	\right.}
\end{equation}
\end{lemma}
In particular, taking $P=P_0$ or $P=\hat P_n$ in Lemma~1, we can find two polytopes $P^*$ and $\tilde P_n$ in $\mathcal P_{r}^{(n)}$ such that
$${ \left\{
    \begin{array}{l}
        |P^*\triangle P_0|\leq \frac{K_1}{n}     \vspace{3mm}\\
        P^*\subseteq P_0^{\sqrt d/n}, \mbox{ }\mbox{ }P_0\subseteq (P^*)^{\sqrt d/n}
    \end{array}
	\right.}
$$
and
$${ \left\{
    \begin{array}{l}
       |\tilde P_n\triangle \hat P_n|\leq \frac{K_1}{n}     \vspace{3mm}\\
        \tilde P_n\subseteq \hat P_n^{\sqrt d/n}, \mbox{ }\mbox{ }\hat P_n\subseteq \tilde P_n^{\sqrt d/n}.
    \end{array}
	\right.}
$$
Note that $\tilde P_n$ is random. Let $\epsilon>0$.
By construction, $|\hat P|\leq |P_0|$, so $|\hat P_n\triangle P_0|\leq 2 |P_0\backslash \hat P_n|$. Besides, if $G_1$, $G_2$ and $G_3$ are three measurable subsets of $\R^d$, the following triangle inequality holds :
\begin{equation}
	\label{triangleineq}|G_1\backslash G_3|\leq |G_1\backslash G_2|+|G_2\backslash G_3|.
\end{equation}
Let us now write the following inclusions between the events.
\begin{align}
	\left\{|\hat P_n\triangle P_0|>\epsilon \right\} & \subseteq \left\{|P_0\backslash \hat P_n|>\epsilon/2 \right\} \nonumber \\
	& \subseteq  \left\{|P^*\backslash \tilde P_n|>\epsilon/2 - \frac{2K_1}{n} \right\} \nonumber \\
	\label{proof1}& \subseteq \bigcup_P  \left\{\tilde P_n = P \right\},
\end{align}
where the latest union is over the class of all $P\in\mathcal P_{r}^{(n)}$ that satisfy the inequality $|P^*\backslash P|>\epsilon/2 - \frac{2K_1}{n}$. Let $P$ be such a polytope, then if $\tilde P_n = P$, then necessarily the sample $\{X_1,\ldots,X_n\}$ is included in $P^{\frac{\sqrt d}{n}}$, by definition of $\tilde P_n$, and \eqref{proof1} becomes
\begin{align}
	\PP_{P_0}\left[\tilde P_n = P \right] & \leq \PP_{P_0}\left[X_i \in P^{\frac{\sqrt d}{n}}, i=1,\ldots,n \right] \nonumber \\
	& \leq \left (1-\frac{|P_0\backslash P^{\frac{\sqrt d}{n}}|}{|P_0|}\right)^n \nonumber \\
	& \leq \left (1-|P_0\backslash P^{\frac{\sqrt d}{n}}|\right)^n, \mbox{ }\mbox{ since } |P_0|\leq 1 \nonumber \\
	& \leq \left (1-|P^*\backslash P|+|P^*\backslash P_0|+|P^{\frac{\sqrt d}{n}}\backslash P|\right)^n, \mbox{ }\mbox{ using }\eqref{triangleineq} \nonumber \\
	& \leq \left (1-\epsilon/2 + \frac{4K_1}{n} \right)^n \nonumber \\
	\label{proof2}& \leq C_1\exp(-n\epsilon/2),
\end{align}
where  $C_1=\mathrm e^{4K_1}$. Therefore, using \eqref{proof1} and \eqref{proof2} and denoting by $\# \mathcal P_{r}^{(n)}$ the cardinality of the finite class $\mathcal P_{r}^{(n)}$,
\begin{align}
	P_{P_0}\left[|\hat P_n\triangle P_0|>\epsilon \right] & \leq \# \mathcal P_{r}^{(n)}C_1\exp(-n\epsilon/2) \nonumber \\
	& \leq (n+1)^{dr}C_1\exp(-n\epsilon/2) \nonumber \\
	\label{proof3}& \leq C_1\exp(-n\epsilon/2+2dr\ln n).
\end{align}
It turns out that if we take $\epsilon$ of the form $\displaystyle{\frac{4dr\ln n}{n}+\frac{x}{n}}$, \eqref{proof3} becomes
\begin{equation}
	\mathbb P_{P_0}\left[n\left(|\hat P_n\triangle P_0|-\frac{4dr\ln n}{n}\right)\geq x\right] \leq C_1\mathrm e^{-x/2},
\end{equation}
which holds for any $x>0$ and any $P_0\in\mathcal P_r$. Theorem \ref{Theorem1} is proved.

Corollary \ref{Corollary1} comes by applying Fubini's theorem (see \cite{Brunel13} for details).

\subsection{Proof of Theorem \ref{theorem5}}

Let $r\geq 10$ be an integer, supposed to be even without loss of generality and assume $n\geq r$. Consider a regular convex polytope $P^*$ in $[0,1]^2$ with center $C=(1/2,1/2)$ and with $r/2$ vertices, denoted by $A_0, A_2, \ldots, A_{r-2}$, such that for all $k=0,\ldots,r/2-1$, the distance between $A_{2k}$ and the center $C$ is $1/2$. Let $A_1, A_3,\ldots, A_{r-1}$ be $r/2$ points built as in Figure~1: for $k=0,\ldots, r/2-1, A_{2k+1}$ is on the mediator of the segment $[A_{2k},A_{2k+2}]$, outside $P^*$, at a distance $\delta=h/2\cos(2\pi/r)\tan(4\pi/r)$ of $P^*$, with $h\in (0,1)$ to be chosen. Note that by our construction, $A_{2k}$ and $A_{2k+2}$ are vertices of the convex hull of $A_0,A_2,\ldots,A_{r-2}$ and $A_{2k+1}$.

\begin{figure}
\centering
\includegraphics[width=7cm]{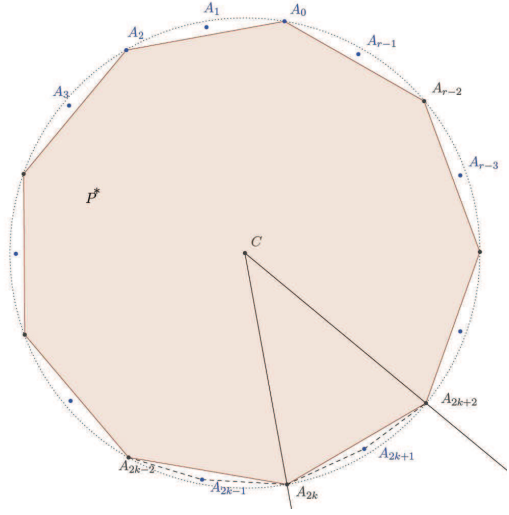}
\caption{Construction of hypotheses for the lower bound}
\label{Fig 1}
\end{figure}

Let us denote by $D_k$ the smallest convex cone with apex $C$, containing the points $A_{2k}, A_{2k+1}$ and $A_{2k+2}$, as drawn in Figure~1. For $\omega=(\omega_0,\ldots,\omega_{r/2-1})\in\{0,1\}^{r/2}$, we denote by $P_\omega$ the convex hull of $P^*$ and the points $A_{2k+1}, k=0,\ldots,r/2-1$ such that $\omega_k=1$. Then we follow the scheme of the proof of Theorem 5 in \cite{Brunel13}.

For $k=0,\ldots,R/2-1$, and $(\omega_0,\ldots,\omega_{k-1},\omega_{k+1},\ldots,\omega_{r/2-1})\in\{0,1\}^{r/2-1}$, we denote by
\begin{align*}
\omega^{(k,0)} &=(\omega_1,\ldots,\omega_{k-1},0,\omega_{k+1},\ldots,\omega_{r/2-1}) \mbox{  and by} \\
\omega^{(k,1)}&=(\omega_1,\ldots,\omega_{k-1},1,\omega_{k+1},\ldots,\omega_{r/2-1}).
\end{align*}

Note that for $k=0,\ldots,r/2-1$, and $(\omega_1,\ldots,\omega_{k-1},\omega_{k+1},\ldots,\omega_{r/2-1})\in$ \linebreak $\{0,1\}^{r/2-1}$,
\begin{equation*}
	|P_{\omega^{(k,0)}}\triangle P_{\omega^{(k,1)}}|=\frac{\delta}{2}\cos(2\pi/r).
\end{equation*}

Let $H$ be the Hellinger distance between probability measures. For the definition and some properties, see \cite{30}, Section 2.4. We have, by a simple computation,
\begin{align*}
	1-\frac{H(P_{\omega^{(k,0)}},P_{\omega^{(j,1)}})^2}{2} & = \sqrt{1-\frac{|P_{\omega^{(k,1)}}\backslash P_{\omega^{(k,0)}}|}{|P_{\omega^{(k,1)}}|}} \nonumber \\
	& = \sqrt{1-\frac{\delta/2\cos(2\pi/r)}{|P_{\omega^{(k,1)}}|}} \nonumber \\
	& \geq \sqrt{1-\frac{\delta\cos(2\pi/r)}{4}}
\end{align*}
since $|P_{\omega^{(k,1)}}|\geq |P^*|\geq 1/2$.
Now, let $\hat P_n$ be any estimator of $P^*$, based on a sample of $n$ i.i.d. random variables. By the same computation as in the proof of Theorem 5 in \cite{Brunel13}, based on the $2^{r/2}$ hypotheses that we constructed, we get

\begin{align}
	\sup_{P\in\mathcal P_r}\mathbb E_P\left[|P\triangle \hat P_n|\right] & \geq \frac{r\delta\cos(2\pi/r)}{8}\left(1-\frac{\delta\cos(2\pi/r)}{4}\right)^n \nonumber \\
	\label{Der1} & \geq \frac{rh\cos\left(\frac{2\pi}{r}\right)^2\tan\left(\frac{4\pi}{r}\right)}{8}\left(1-\frac{h\cos\left(\frac{2\pi}{r}\right)^2\tan\left(\frac{4\pi}{r}\right)}{8}\right)^n.
\end{align}

Note that if we denote by $\displaystyle{x=\frac{2\pi}{r}>0}$ and $\displaystyle{\phi(x)=\frac{1}{x}\cos(x)\tan(2x)}$, then $\phi(x)\gtrsim 1$ since $r$ is supposed to be greater or equal to $10$. Therefore, by the choice $h=r/n\leq 1$ (we assumed that $n\geq r$), \eqref{Der1} becomes

\begin{equation*}
	\sup_{P\in\mathcal P_r}\mathbb E_P\left[|P\triangle \hat P_n|\right] \gtrsim \frac{r}{n},
\end{equation*}
and Theorem \ref{theorem5} is proved.

\subsection{Proof of Theorem \ref{theorem3}}

Let $t>0$ be fixed. Let $G_1=t^{1/d}B_2^d$ and $G_2=(2t)^{1/d}B_2^d$.
Let us denote respectively by $\PP_1$ and $\PP_2$ the uniform distributions on $G_1$ and $G_2$, and by $\E_1$ and $\E_2$ the corresponding expectations. We denote by $\PP_1^{\otimes n}$ and $\PP_2^{\otimes n}$ the $n$-product of $\PP_1$ and $\PP_2$, respectively, i.e. the probability distribution of a sample of $n$ i.i.d. random variables of distribution $\PP_1$ and $\PP_2$, respectively. The corresponding expectations are still denoted by $\E_1$ and $\E_2$.
Then, for any estimator $\hat G_n$ based on a sample of $n$ random variables, we bound from bellow the minimax risk by the Bayesian one.
\begin{align*}
	\sup_{G\in\mathcal K}\E_G[|\hat G_n\triangle G|] & \geq \frac{1}{2}\left(\E_1[|\hat G_n\triangle G_1|]+\E_2[|\hat G_n\triangle G_2|]\right) \\
	& \geq \frac{1}{2}\int_{(\R^d)^n}\left(|\hat G_n\triangle G_1|+|\hat G_n\triangle G_2|\right)\min(d\PP_1^{\otimes n},d\PP_2^{\otimes n}) \\
	& \geq \frac{1}{2}\int_{(\R^d)^n}|G_1\triangle G_2|\min(d\PP_1^{\otimes n},d\PP_2^{\otimes n}) \\
	& \geq \frac{t|B_2^d|}{4}\left(1-\frac{H(\PP_1^{\otimes n},\PP_2^{\otimes n})^2}{2}\right)^2,
\end{align*}
where $H$ is the Hellinger distance between probability measures, as in the proof of Theorem \ref{theorem5}.
Therefore,
\begin{equation}
	\label{proof3step1}\sup_{G\in\mathcal K}\E_G[|\hat G_n\triangle G|] \geq \frac{t|B_2^d|}{4}\left(1-\frac{H(\PP_1,\PP_2)^2}{2}\right)^{2n},
\end{equation}
and a simple computation shows that
\begin{equation*}
	1-\frac{H(\PP_1,\PP_2)^2}{2} = \frac{1}{\sqrt{2}},
\end{equation*}
and \eqref{proof3step1} becomes
\begin{equation*}
	\sup_{G\in\mathcal K}\E_G[|\hat G_n\triangle G|] \geq \frac{t|B_2^d|}{2^{n+2}}.
\end{equation*}
This ends the proof of Theorem \ref{theorem3} by taking $t$ arbitrarily large.

\subsection{Proof of Theorem \ref{Theorem2}}

We first state the following result, due to Groemer \cite{15bis}.
\begin{lemma}
Let $K\in\mathcal K$ and $n$ be an integer greater than the dimension $d$. Let $\hat K_n$ be the convex hull of a sample of $n$ i.i.d. random variables uniformly distributed in $K$. Then, if $B$ denotes a Euclidean ball in $\R^d$, of the same volume as $K$,
\begin{equation*}
	\E_K[|K\triangle\hat K_n|]\leq \E_B[|B\triangle\hat K_n|].
\end{equation*}
\end{lemma}
If $K\in\mathcal K$, we denote by $K_1=\frac{1}{|K|^{1/d}}K$, so $K_1$ is homothetic to $K$ and has volume 1. Besides, if $X_1,\ldots,X_n$ are i.i.d. with uniform distribution in $K$, then $\frac{1}{|K|^{1/d}}X_1,\ldots,\frac{1}{|K|^{1/d}}X_n$ are i.i.d. with uniform distribution in $K_1$, and
\begin{equation*}
	\E_K\left[\frac{|K\triangle\hat K_n|}{|K|}\right] = \E_{K_1}[|K_1\triangle\hat K_n|]
\end{equation*}
Therefore, one gets, from the previous lemma, that if $B=B_2^d$ is the unit Euclidean ball in $\R^d$,
\begin{equation*}
	\E_K\left[\frac{|K\triangle\hat K_n|}{|K|}\right]\leq \E_B\left[\frac{|B\triangle\hat K_n|}{|B|}\right], \forall K\in\mathcal K.
\end{equation*}
In other terms,
\begin{equation*}
	\mathcal Q_n(\tilde K_n ; \mathcal K)\leq \E_B\left[\frac{|B\triangle\hat K_n|}{|B|}\right].
\end{equation*}
By \cite{3'''0}, $\displaystyle{n^\frac{2}{d+1}\E_B\left[\frac{|B\triangle\hat K_n|}{|B|}\right]}$ tends to the affine surface area of $B=B_2^d$ times some positive constant which depends on $d$ only, that is,
\begin{equation*}
	n^\frac{2}{d+1}\E_B\left[\frac{|B\triangle\hat K_n|}{|B|}\right] \rightarrow C(d),
\end{equation*}
where $C(d)$ is a positive constant which depends on $d$ only. Thus, Theorem \ref{Theorem2} is proved.

\subsection{Proof of Theorem \ref{theorem4}}

Let $\hat r$ be chosen as in Section \ref{sec:4}.
Let $d+1\leq r^*\leq R_n$ and $P\in\mathcal P_{r^*}$. We distinguish two cases,
\begin{equation}
	\label{p_adapt2} \mathbb E_P[|\hat P_n^{adapt}\triangle P|]  =\mathbb E_P[|\hat P_n^{adapt}\triangle P|I(\hat r\leq r^*)] +  \mathbb E_P[|\hat P_n^{adapt}\triangle P|I(\hat r> r^*)],
\end{equation}
and we bound separately the two terms in the right side. Note that if $\hat r\leq r^*$, then, $\hat r\leq R_n$, so $\hat P_n^{adapt}=\hat P_n^{(\hat r)}$, and by definition of $\hat r$
\begin{equation}
	\label{adapt111} |\hat P_n^{(r^*)}\triangle\hat P_n^{(\hat r)}|\leq \frac{Cr^*\ln n}{n}.
\end{equation}
Therefore, using the triangle inequality,
\begin{align}
	\mathbb E_P &[|\hat P_n^{adapt}\triangle P|I(\hat r\leq r^*)] \nonumber \\
	& \leq \mathbb E_P[|\hat P_n^{adapt}\triangle \hat P_n^{(r^*)}|I(\hat r\leq r^*)] + \mathbb E_P[|\hat P_n^{(r^*)}\triangle P|I(\hat r\leq r^*)] \nonumber \\
	\label{p_adapt5}& \lesssim_ d \frac{r^*\ln n}{n}, \mbox{ }\mbox{by \eqref{adapt111} and Corollary~1}.
\end{align}
The second term of \eqref{p_adapt2} is bounded differently. First note that $\hat P_n^{adapt}\subseteq[0,1]^d$, so $P\triangle\hat P_n^{adapt}\subseteq[0,1]^d$ and $|P\triangle\hat P_n^{adapt}|\leq 1$ almost surely. Besides, note that if $\hat r> r^*$, then for some $r\in\{r^*+1,\ldots,n\}, |\hat P_n^{(r^*)}\triangle\hat P_n^{(r)}|> \frac{Cr\ln n}{2n}$. Otherwise, for any $r_1, r_2\in\{r^*,\ldots,n\}$, one would have, by the triangle inequality, $ |\hat P_n^{(r_1)}\triangle\hat P_n^{(r_2)}|\leq\frac{Cr\ln n}{n}$, and this would contradict that $\hat r > r^*$. Thus, we have the following inequalities.

\begin{align}
	 \mathbb E_P[|\hat P_n^{adapt}\triangle P|I(\hat r> r^*)] & \leq \mathbb P_P[\hat r> r^*] \nonumber\\
	 & \leq \sum_{r=r^*+1}^n\mathbb P_P\left[|\hat P_n^{(r^*)}\triangle\hat P_n^{(r)}|> \frac{Cr\ln n}{2n}\right] \nonumber
\end{align}
\begin{align}
	 & \leq \sum_{r=r^*+1}^n\mathbb P_P\left[|\hat P_n^{(r^*)}\triangle P|+|\hat P_n^{(r)}\triangle P|> \frac{Cr\ln n}{2n}\right] \nonumber\\
	 & \label{p_adapt3} \leq \sum_{r=r^*+1}^n\left (\mathbb P_P\left[|\hat P_n^{(r^*)}\triangle P|> \frac{Cr\ln n}{4n}\right]+\mathbb P_P\left[|\hat P_n^{(r)}\triangle P|> \frac{Cr\ln n}{4n}\right]\right).
\end{align}

Note that since $P\in\mathcal P_{r^*}$, it is also true that $P\in\mathcal P_{r}, \forall r\geq r^*$. Therefore, for $d+1\leq r^*\leq r\leq n$, we have, using Theorem~1, with $x=(C/4-4d)r\ln n$,
\begin{equation*}
    \mathbb P_P\left[|\hat P_n^{(r)}\triangle P|> \frac{Cr\ln n}{4n}\right] \leq e^{-(C/8-2d)r\ln n}\leq n^{-(C/8-2d)(d+1)}\lesssim_ d n^{-1},
\end{equation*}
by the choice of $C$.

It comes from \eqref{p_adapt3} that
\begin{equation}
	\label{p_adapt4} \mathbb E_P[|\hat P_n^{adapt}\triangle P|I(\hat r> r^*)] \lesssim_ d n^{-1}.
\end{equation}

Finally, using \eqref{p_adapt5} and \eqref{p_adapt4},
\begin{equation*}
	\mathbb E_P[|\hat P_n^{adapt}\triangle P|] \lesssim_ d \frac{r^*\ln n}{n}.
\end{equation*}

Let us now assume that the unknown support, which we now denote by $K$, is any convex body in $\mathcal P_\infty=\mathcal K_1$, possibly a polytope with many (more than $R_n$) vertices.
We write, similarly to the previous case,
\begin{equation}
	\label{k_adapt2} \mathbb E_K[|\hat P_n^{adapt}\triangle K|]  = \mathbb E_K[|\hat P_n^{adapt}\triangle K|I(\hat r\leq R_n)] +  \mathbb E_K[|\hat P_n^{adapt}\triangle K|I(\hat r> R_n)],
\end{equation}
and we bound separately the two terms of the right side.
If $\hat r\leq R_n$, then $\hat P_n^{adapt}=\hat P_n^{(\hat r)}$. As we already explained, if $V_n$ is the number of vertices of $\hat K_n$, then $\hat r\leq V_n$ almost surely, and $\hat K_n=\hat P_n^{(V_n)}$. So we have, $|\hat P_n^{adapt}\triangle\hat K_n|=|\hat P_n^{(\hat r)}\triangle\hat P_n^{(V_n)}|\leq \frac{CV_n\ln n}{n}$. Therefore, using the triangle inequality,
\begin{align}
	\mathbb E_K[|\hat P_n^{adapt} & \triangle K|I(\hat r \leq R_n)] \leq \mathbb E_K[|\hat P_n^{adapt}\triangle \hat K_n|I(\hat r \leq R_n)] \nonumber \\
	& \mbox{ }\mbox{ }\mbox{ }\mbox{ }\mbox{ }\mbox{ }\mbox{ }\mbox{ }\mbox{ }\mbox{ }\mbox{ }\mbox{ }\mbox{ }\mbox{ }\mbox{ }\mbox{ }\mbox{ }\mbox{ }\mbox{ }\mbox{ }\mbox{ }\mbox{ }\mbox{ }\mbox{ }+\mathbb E_K[|\hat K_n\triangle K|I(\hat r \leq R_n)]  \nonumber \\
	& \leq \mathbb E_K\left[\frac{CV_n\ln n}{n}I(\hat r \leq R_n)\right]+\mathbb E_K[|\hat K_n\triangle K|] \nonumber \\
	\label{k_adapt3} & \leq \frac{C\mathbb E_K[V_n]\ln n}{n}+\mathbb E_K[|\hat K_n\triangle K|].
\end{align}

We have the following lemma, which is a consequence of Efron's equality \eqref{EfronId} and Theorem \ref{Theorem2}.

\begin{lemma}
	$$\sup_{K\in\mathcal K}\E_K[V_n]\lesssim_ d n^\frac{d-1}{d+1}.$$
\end{lemma}

Therefore, \eqref{k_adapt3} becomes, using Theorem \ref{Theorem2},
\begin{equation}
	\label{k_adapt4} \mathbb E_K[|\hat P_n^{adapt}\triangle K|I(\hat r \leq R_n)] \lesssim_ d (\ln n)n^{-\frac{2}{d+1}}.
\end{equation}

The second term is easily bounded. If $\hat r> R_n$, then $\hat P_n^{adapt}=\hat K_n$ and
\begin{align}
	\mathbb E_K[|\hat P_n^{adapt}\triangle K|I(\hat r > R_n)] & \leq  \mathbb E_K[|\hat K_n\triangle K|] \nonumber \\
	\label{k_adapt5} & \lesssim_ d n^{-\frac{2}{d+1}},
\end{align}
by Theorem \ref{Theorem2}.

From \eqref{k_adapt4} and \eqref{k_adapt5} we get
\begin{equation*}
	\mathbb E_K[|\hat P_n^{adapt}\triangle K|] \lesssim_ d (\ln n)n^{-\frac{2}{d+1}}.
\end{equation*}

Theorem \ref{theorem4} is then proven.



\end{document}